\documentclass[1p]{elsarticle}

\usepackage{graphicx}
\usepackage{comment}
\usepackage{amssymb}
\usepackage{amsmath}
\usepackage{latexsym}
\usepackage[utf8]{inputenc}
\usepackage{algorithm}
\usepackage{algpseudocode}
\usepackage{float}
\usepackage[normalem]{ulem} 
\usepackage[caption = false]{subfig}
\usepackage{graphicx}
\usepackage{commath}
\usepackage{caption}
\usepackage{amsthm}
\usepackage{url}
\DeclareMathOperator{\E}{\mathbb{E}}

\newtheorem{theorem}{Theorem}
\makeatletter
\def\ps@pprintTitle{%
 \let\@oddhead\@empty
 \let\@evenhead\@empty
 \def\@oddfoot{}%
 \let\@evenfoot\@oddfoot}
\makeatother

\begin{document}
\begin{frontmatter}

\title{Error estimators and their analysis for CG, Bi-CG and GMRES}

%% or include affiliations in footnotes:\cortext[mycorrespondingauthor]{Corresponding author}
\author[]{Puneet Jain}
\author[]{Krishna Manglani}

\author[]{Murugesan Venkatapathi}
\ead{murugesh@iisc.ac.in}

\address{Department of Computational \& Data Sciences, Indian Institute of Science, Bangalore - 560012}

\begin{abstract}
The demands of accuracy in measurements and engineering models today, renders the condition number of problems larger. While a corresponding increase in the precision of floating point numbers ensured a stable computing, the uncertainty in convergence when using residue as a stopping criterion has increased. We present an analysis of the uncertainty in convergence when using relative residue as a stopping criterion for iterative solution of linear systems, and the resulting over/under computation for a given tolerance in error. This shows that error estimation is significant for an efficient or accurate solution even when the condition number of the matrix is not large. An $\mathcal{O}(1)$ error estimator for iterations of the CG algorithm was proposed more than two decades ago. Recently, an $\mathcal{O}(k^2)$ error  estimator was described for the GMRES algorithm which allows for non-symmetric linear systems as well, where $k$ is the iteration number. We suggest a minor modification in this GMRES error estimation for increased stability. In this work, we also propose an $\mathcal{O}(n)$ error estimator for A-norm and $l_{2}$ norm of the error vector in Bi-CG algorithm. The robust performance of these estimates as a stopping criterion results in increased savings and accuracy in computation, as condition number and size of problems increase.
\end{abstract}


\end{frontmatter}

\section{Introduction}

Solving a system of linear equations in the form $Ax = b$ is a ubiquitous requirement in science and engineering (where $A$ is a given matrix, $x$ and $b$ are the unknown and known vectors respectively; $x$ $\in \mathbb{R}^n$ and $b$ $\in \mathbb{R}^n$ if $A$ $\in \mathbb{R}^{n \times n}$). Iterative methods like Conjugate Gradient (CG) \cite{hestenes1952methods}, Bi-Conjugate Gradient (Bi-CG) \cite{fletcher1976conjugate}, and Generalized Minimal Residual (GMRES) \cite{saad1986gmres}, mostly known as Krylov subspace methods, are commonly used to solve large linear systems.
While using such methods, iterations are to be stopped when the norm of the error ($\epsilon_{k} = x - x_{k}$) is less than a desired tolerance, where $x$ is the unknown solution to the linear system and $x_k$ is its iterative approximation at the $k^{th}$ iteration. Since the actual error is unknown, relative residue ($\frac{\lVert r_k \rVert}{\lVert b \rVert}$) is considered as a stopping criteria where $r_{k} = b - Ax_{k}$ is the residual vector at $k^{th}$ iteration. Such stopping criteria are useful when the system is well-conditioned but they can become unreliable for even a moderate increase in the relative condition number of $A$ ($\kappa \gtrsim 10$, where $\kappa = \norm{A}_{2}\norm{A^{-1}}_{2}$), depending on the choice of initial approximation. It is well known that relying on relative residue may lead one to stop the iterations too early when the norm of error is still much larger than tolerance, or not stop early enough with too many floating point operations done for the required accuracy (\cite{meurant2013computing}; \cite{hestenes1952methods}). Moreover, when condition number of the matrix is large, the residual of algorithms like Bi-CG need not show monotonic behaviour and oscillate \cite{bjorck2015numerical} while the actual error might still be (however slowly) converging, and similarly, oscillations of errors in the GMRES algorithm are possible when the residue converges slowly. Note that the $l_{2}$ norm of the relative errors and residues can differ significantly as given below:

\begin{equation}
   \kappa \frac{\lVert r_k \rVert}{\lVert b \rVert} \geq \frac{\lVert \epsilon_k \rVert}{\lVert x \rVert} \geq \frac{\lVert r_k \rVert}{\lVert b \rVert} \frac{1}{\kappa}
\end{equation}

Even when the iterative algorithms are used with preconditioners that reduce the condition number of the matrix, a reduction of the condition number of the backward problem to be solved is not guaranteed, as required for the faster convergence (see section \ref{sec:Cond_numbers}). This is observed with matrices  of larger dimensions or higher condition numbers, making the preconditioners ineffective in many cases. Moreover, most of the preconditioning techniques often fail when matrices are non-symmetric and indefinite \cite{saad1988preconditioning}. Thus, for even $\kappa(A) \gtrsim 10^2$, either the accuracy or the efficiency of computation is degraded by the above conundrum. The requirements of precision in measurements and engineering today renders both the size and condition number of most problems large; making accurate stopping and restarting criteria indispensable in ensuring computational efficiency of solvers. This motivated methods to estimate norms of the error in iterative solvers. Approximations of the errors in $x_k$ for an iterative method were suggested in general more than two decades ago \cite{brezinski1999error} \cite{brezinski2008error}. Error estimators for the CG algorithm were also presented \cite{meurant2006lanczos}, and for solving non-symmetric linear systems using Full-Orthogonalization (FOM) and GMRES methods, formulas for estimation of errors have been suggested more recently \cite{meurant2011estimates}. We suggest a minor modification to this GMRES estimator to increase its stability and precision. We also derive an efficient estimator of error in solving non-symmetric linear systems using Bi-CG.

In the second part of this work, we present an analysis highlighting the significance of these error estimators for CG, Bi-CG and GMRES methods. We confirm that there is indeed less certainty in relying on relative residual as an estimator of the true relative error, in comparison to using an error estimator. We quantify this `uncertainty' and show that its \textit{expected} value is proportional to condition numbers i.e. forward or backward condition numbers or their weighted averages depending on the algorithm used. Numerical results averaged over a large number of linear systems were used to demonstrate the theorem on this uncertainty in convergence. The under and over computation due to the uncertainty is demonstrated, and the gains due to the estimators are relatively straight-forward to infer.  Our analysis shows that these error estimators are robust and increase the efficiency/accuracy of computing notably. Section \ref{section2} briefly introduces related work and the previous research in the estimation of error in CG and GMRES. Section \ref{Section3} presents the modification to the GMRES estimator, the proposed estimates of error in Bi-CG, and the analysis of error estimates. Section \ref{Section4} discusses the experimental results on the modified GMRES error estimates, proposed BiCG error estimates, and the performance of error estimates as an efficient stopping criteria for these iterative methods of solving linear systems.

\section{Related work} \label{section2}

\subsection{Error estimators for CG algorithm}

One of the most commonly used methods for solving linear systems with a real symmetric positive definite matrix is the Conjugate Gradient (CG) algorithm. CGQL (CG-Qaudrature Lanczos) algorithm \cite{golub1997matrices, meurant2006lanczos, meurant2013computing} provides bounds of A-norm of error i.e $\norm{x-x_k}_{A}$ at each iteration of CG by using Lanczos and quadrature rules with a delay of $d$ iterations in the estimation. Bounds to the square of the A-norm of error at iteration $k$ of CG are given by \cite{golub1997matrices}:

\begin{equation}\label{lanc_tri}
\lVert \epsilon_{k} \rVert_{A}^{2} \approx \lVert r_{0} \rVert ^{2}_2 [(T_{k+d}^{-1})_{(1,1)} - (T_{k}^{-1})_{(1,1)} ]
\end{equation}

Here, $T_{k}$ is the tridiagonal matrix produced by the Lanczos algorithm at the $k^{th}$ iteration whose entries can be computed from the coefficients of CG algorithm equally well. $(T_k^{-1})_{(1,1)}$ represents the $(1,1)$ element of the matrix $T_{k}^{-1}$. Also $\epsilon$ and $r$ are the error and residual vectors of the CG algorithm respectively. The essence of CGQL lies in computing difference between $(1,1)$ elements of the inverse of two tridiagonal matrices generated from a Lanczos algorithm with the same starting vectors as the CG algorithm. However, from a computational point of view, other equivalent formulas for approximation of error in CG were proposed \cite{meurant1997computation}, including a $LDL^{T}$ factorization of the corresponding Jacobi matrix available in CG \cite{meurant2013computing}. Estimates of the $l_2$ norm of the error in the CG algorithm were later proposed \cite{meurant2005estimates}, including extensions to preconditioned versions \cite{strakovs2005error}. It should also be noted that with a delay in estimation $d$, where $d>1$, considering $\norm{\epsilon_{k+d}}_A \ll \norm{\epsilon_{k}}_A$, the approximation of the A-norm of error at iteration $k$ is also given using CG coefficients $\alpha_j$, as: 

\begin{equation}
\nonumber
\lVert \epsilon_{k} \rVert_{A}^{2} \approx \sum_{j=k}^{k+d-1} \alpha_{j} r_{j}^{T} r_{j}
\end{equation}

\subsection{Error estimator for GMRES algorithm}

Generalized minimal residual (GMRES) algorithm is an iterative method for solving linear systems even for a general non-symmetric matrix \cite{saad1986gmres}. Let $V_k$ be a matrix whose columns are orthonormal basis vectors $v_j,\ j=1 \dots k$, of Krylov subspace $\mathcal{K}_k(A,r_0)$, where $r_0$ is initial residual. The iterates of GMRES are defined as  $x_k = x_0 + V_k z_k$ where vector $z_k$ has weights for the orthonormal basis vectors $v_j$ in the matrix $V_k$. Maintaining the notations \cite{meurant2011estimates} for estimates of norm of error in FOM and GMRES at the $k^{th}$ iteration of GMRES, the Hessenberg matrix $H_k$ resulting from the iterative orthogonal transformation of the system matrix can be decomposed block-wise as:

\[
H_k = 
\begin{pmatrix}
H_{k-d} & W_{k-d} \\
Y_{k-d}^T & \tilde{H}_{k-d}
\end{pmatrix}
\]

Here, the matrices $H_{k-d}$ and $\tilde{H}_{k-d}$ are square and upper Hessenberg and $W_{k-d}$ is generally a full rectangular matrix. Let $h_{i,j}$ be the $(i,j)$ element of the Hessenberg matrix, $e_k$ is the $k^{th}$ column of an identity matrix, the vector $t_k$ be the last column of $\big( H_k^T H_k  \big)^{-1} $ with $t_{k,k}$ as its last element, then variables $\gamma_{k-d}, w_{k-d},$ $\delta_{k+1}$, and $u_k$ can be defined as the following,
\[
    \gamma_{k-d} = \dfrac{h_{k-d+1,k-d}\big( e_{k-d}, H_{k-d}^{-1}e_1 \big)}{1 - h_{k-d+1,k-d}\big( e_{k-d}, H_{k-d}^{-1}w_{k-d} \big)}
\]
\[
    w_{k-d} = W_{k-d}\tilde{H}_{k-d}^{-1}e_1
\]
\[
    \delta_{k+1} = \dfrac{h_{k+1,k}^2}{1+h_{k+1,k}^2 t_{kk}}
\]

\[
    u_k = \delta_{k+1}t_k
\]

The error estimate $\epsilon_{k-d}$ at $(k-d)^{th}$ iteration of GMRES can be aproximated by \cite{meurant2011estimates}:

\begin{equation} \label{eq:Estimator_Meurant}
\dfrac{\norm{\epsilon_{k-d}}^2}{\norm{r_0}^2} \approx \gamma_{k-d}^2 \norm{\tilde{H}_{k-d}^{-1}e_1}^2  + \norm{\gamma_{k-d} H_{k-d}^{-1}w_{k-d} + (e_{k-d},H_{k-d}^{-1}e_1)u_{k-d}}^2
\end{equation}

Note that the above estimation involves only one dense matrix-vector multiplication and requires $\mathcal{O}(k^2)$ arithmetic operations every iteration, and this is negligible for all $k \ll n$. For non-converging problems, as $k \longrightarrow n$, the total additional cost due to the estimator may approach a third of the cost of the GMRES iterations, but this is likely to happen only rarely in applications with a large $n$, and thus the estimator can provide significant benefits on the average.

\section{Methods}\label{Section3}

\subsection{Proposed modification to the GMRES error estimator} 

The above estimate of error in GMRES can also be rewritten as term-1 in the following:

\begin{equation} \label{estanlysis}
\dfrac{\norm{\epsilon_{k-d}}^2}{\norm{r_0}^2} = \underbrace{\dfrac{\norm{\epsilon_{k-d}}^2 - \norm{\epsilon_{k}}^2}{\norm{r_0}^2}}_\text{term1}  +   \underbrace{\norm{s_k}^2}_\text{term2} +  \underbrace{\zeta_k}_\text{term3}
\end{equation}

where $s_k = ( e_k,H_k^{-1}e_1) u_k $ and\\
\begin{equation*}
\begin{split}
 \zeta_k =  & 2h_{k+1,k} \left( e_k,H_k^{-1}e_1 \right) \left((H_n^{-1}e_{k+1})^k,H_k^{-1}e^1 + s_k \right) \\
& - 2h_{k-d+1,k-d} \left( e_{k-d},H_{k-d}^{-1}e_1 \right) \left( (H_n^{-1}e_{k-d+1})^{k-d} - (H_k^{-1}e_{k-d+1})^{k-d} ,H_{k-d}^{-1}e_1 + s_{k-d}  \right) 
\end{split}
\end{equation*}
where superscripts $k$ and $k-d$ have been introduced for the vector components, indicating the appropriate truncated dimensions.

The R.H.S of \eqref{estanlysis} was divided into three terms where the third term given by $\zeta_k$ is assumed to be negligible. Also, $H_n^{-1}$ required in the evaluation of $\zeta_k$ cannot be computed at the $k^{th}$ iteration. Note that for a large $d$ which is a delay in estimation, the second term $\norm{s_k}^2 $ also becomes negligible. However $d$ is kept small in practice ($\approx$ 10), hence in such cases, this second term can significantly contribute to the value of error estimate. We propose to include this second term in the estimates of the error, and here, an absolute operation in \eqref{newesteq} is necessary to preserve the positive estimate of the $l_2$ norm of error using the truncated expression. Including $s_k$ in \eqref{eq:Estimator_Meurant} and considering the absolute values,

\begin{equation} \label{newesteq}
    \dfrac{\norm{\epsilon_{k-d}}^2}{\norm{r_0}^2} \approx   \left\vert \gamma_{k-d}^2 \norm{\tilde{H}_{k-d}^{-1}e_1}^2  +   \norm{\gamma_{k-d} H_{k-d}^{-1}w_{k-d} + (e_{k-d},H_{k-d}^{-1}e_1)u_{k-d}}^2 - \norm{s_k}^2  \right\vert
\end{equation}.

\subsection{Proposed estimator for A-norm and $l_2$ norm of errors in Bi-Conjugate gradient Algorithm (BiCG)} \label{BiCGSection}

Whereas A-norm of error in CG is minimized monotonically with the number of iterations, A-norm of error in BiCG can oscillate and might not even be defined when A is an indefinite matrix. Hence, we consider the absolute values of the bi-linear products of the residues of BiCG as shown below.

\begin{eqnarray}\label{A-norm}
  \lVert \epsilon_{k} \rVert_{A}^2  = \abs{ \epsilon_{k}^{T}A\epsilon_k } = \abs {r_{k}^{T}(A^{T})^{-1}r_{k}} = \abs{ r_{k}^{T}A^{-1}r_{k}}
  \end{eqnarray}
  
Note than when $A$ $\in \mathbb{R}^{N \times N}$ and $r$ $\in \mathbb{R}^N$ ; $r_{k}^{T}(A^{T})^{-1}r_{k}$ is a scalar quantity that is always positive when $A$ is positive definite. Moreover the square of $l_{2}$ norm of error is given by: 

\begin{equation}\label{l2 norm}
   \lVert \epsilon_{k} \rVert_{2}^{2} = \epsilon_{k}^{T}\epsilon_{k} = r_{k}^{T}(A^{T})^{-1}A^{-1}r_{k}  
\end{equation}

If $A=A^{T}$, A-norm and $l_{2}$ norm of error are given by $r_{k}^{T}A^{-1}r_{k}$ and $r_{k}^{T}A^{-2}r_{k}$ respectively. In the following sections, estimates of A-norm and $l_{2}$ norm of error at every iteration of BiCG are derived. BiCG is included as Algorithm \ref{BiCG algorithm} in the appendix for reference.

\subsubsection{Estimation of A-norm of error}

By writing the consecutive difference between A-norm of error at iteration $k$ and $k+1$, the following relation is obtained when A is a non-symmetric matrix (See Appendix for BiCG Algorithm):

\begin{equation}
\nonumber
\begin{split}
 r_{k}^T (A^{T})^{-1} r_{k} - r_{k+1}^T(A^{T})^{-1} r_{k+1} & =  r_{k}^T (A^{T})^{-1} r_{k} - (r_{k} - \alpha_{k}Ap_{k})^{T} (A^{T})^{-1} (r_{k} - \alpha_{k}Ap_{k})\\
 & =  \alpha_{k} p_{k}^T r_{k} + \alpha_{k}r_{k}^{T}(A^{T})^{-1}Ap_{k}  -   \alpha_{k}^{2} p_{k}^T A p_{k}
\end{split}
\end{equation}

In the above equation, the first and third terms of R.H.S can be trivially computed using iterates of the Bi-CG Algorithm. Second term involves computation of $A^{-1}$ hence we further reduce it. As $ r_{k+1}  = r_{k} - \alpha_{k} A p_{k} $ and $Ap_{k} = \dfrac{r_{k} - r_{k+1}}{\alpha_{k}} $, the following relation is obtained:

\begin{equation}\label{eq:11}
\nonumber
\begin{array}{ccl}
 r_{k}^T A^{-1} r_{k} - r_{k+1}^T A^{-1} r_{k+1} = \alpha_{k} r_{k}^T p_{k} + \alpha_{k} r_{k}^T (A^{T})^{-1} (\dfrac{r_{k} - r_{k+1}}{\alpha_{k}})   -  \alpha_{k}^{2} p_{k}^T A p_{k}\\
 \implies r_{k+1}^T A^{-1} r_{k+1} = - \alpha_{k} r_{k}^T p_{k} +   r_{k+1}^T A^{-1} r_{k}   +  \alpha_{k}^{2} p_{k}^T A p_{k}
 \end{array}
\end{equation}

Error ($\epsilon_{k} = x - x_{k}$) is given by $A^{-1}r_{k}$. Also, $\epsilon_k$ can be written as a weighted sum of search directions from $k$ to $n$ given below in \eqref{eq:13}.

\begin{equation}\label{eq:13}
\begin{array}{ccl}
\epsilon_{k} = \sum_{j=k}^{n} \alpha_{j} p_{j} \implies \epsilon_{k} - \epsilon_{k+d} = \sum_{j=k}^{k+d} \alpha_{j} p_{j}\\
\implies  \epsilon_{k} \approx \sum_{j=k}^{k+d} \alpha_{j} p_{j}
\end{array}
\end{equation}

If $\epsilon_{k+d}$ denotes the error at $k+d$ iteration and $\epsilon_{k+d} \ll \epsilon_{k}$ for some $d > 0$, then further terms of the series can be neglected. The error vector in \eqref{eq:13} can provide an estimate of A-norm of error after inducing a delay of $d$ iterations:

\begin{equation}\label{eq:14}
   \lVert \epsilon_{k+1} \rVert_{A}^2 = r_{k+1}^T A^{-1} r_{k+1} \approx - \alpha_{k} r_{k}^T p_{k} +   r_{k+1}^T(\sum_{j=k}^{k+d}\alpha_{j}p_{j})   +  \alpha_{k}^{2} p_{k}^T A p_{k}
\end{equation}

When A is positive definite, the above expression is always positive and thus provides a lower bound for the square of A-norm of error. Equation \eqref{eq:14} is derived directly from Bi-CG and when $A = A^T$, this estimator for $A$-norm is identical to the earlier proposed CGQL $A$-norm estimator (\ref{BCGtoCG}). Since the matrix-vector product $Ap_k$ is known in terms of the residual vectors, the above evaluation costs only $2n$ arithmetic operations, compared to the $\mathcal{O}(n^2)$ arithmetic operations in each Bi-CG iteration. The summation of $d$ terms in the above implies a latency of $d$ iterations to begin with, when an estimate of the error is not available. 

\subsubsection{ Estimation of $l_{2}$ norm of error}

By writing the consecutive difference between $l_{2}$ norm of error at iteration $k$ and $k+1$ of Bi-CG, the following relation is obtained when A is a non-symmetric matrix:

\begin{multline}\label{eq:15}
 r_{k}^T (A^{T})^{-1} A^{-1} r_{k} - r_{k+1}^T(A^{T})^{-1}A^{-1} r_{k+1} =  r_{k}^T (A^{T})^{-1} A^{-1}r_{k}\\ - (r_{k} - \alpha_{k}Ap_{k})^{T} (A^{T})^{-1} A^{-1}(r_{k} - \alpha_{k}Ap_{k})   
\end{multline}

Rearranging \eqref{eq:15},

\begin{equation}\label{eq:20}
     r_{k}^T (A^{T})^{-1} A^{-1} r_{k} - r_{k+1}^T(A^{T})^{-1}A^{-1} r_{k+1} = 2\alpha_{k} p_{k}^T A^{-1}r_{k}   -   \alpha_{k}^{2} \lVert p_{k} \rVert ^{2}
\end{equation}

Using \eqref{eq:13} and $\epsilon_{k} = A^{-1}r_{k}$, the above equation can be rewritten as:

\begin{equation}\label{eq:18}
 \lVert \epsilon_{k+1} \rVert_{2}^2 = r_{k+1}^T(A^{T})^{-1}A^{-1} r_{k+1}  \approx (\sum_{j=k}^{k+d} \alpha_{j} p_{j})^{T}(\sum_{j=k}^{k+d} \alpha_{j} p_{j}) -2\alpha_{k} p_{k}^T (\sum_{j=k}^{k+d} \alpha_{j} p_{j}) +   \alpha_{k}^{2} \lVert p_{k} \rVert ^{2} 
\end{equation}

The evaluation involves only $2n$ arithmetic operations every iteration, and for the first $d$ iterations an estimate of the error is not available. This estimation adds a negligible computational overhead compared to the $\mathcal{O}(n^2)$ cost of every Bi-CG iteration, especially as $n$ increases. 

\subsection{Analysis of estimators}

\subsubsection{Condition number of the problem} \label{sec:Cond_numbers}

Condition numbers of problems enable us to estimate the expected accuracy of computed result. In this case, the relative condition number for computing $b$ given $A$ and $x$ (forward problem), and condition number for computing $x$ from $A$ and $b$ (backward problem) are given by \eqref{cn of prob1} and \eqref{cn of prob2} respectively \cite{trefethen1997numerical}. They help us infer the \textit{relative} change in the output for a relative change in inputs of the problem representing round-offs, thus constraining the accuracy of the output. The norm ($\norm{.}$), unless specified, is considered to be the $l_2$ norm.

\begin{equation}\label{cn of prob1}
    \kappa_{f}(A,x) = \lVert A \rVert \frac{\lVert x \rVert}{\lVert Ax \rVert}
\end{equation}

\begin{equation}\label{cn of prob2}
    \kappa_{b}(A,b) = \lVert A^{-1} \rVert \dfrac{\lVert b \rVert}{\lVert A^{-1}b \rVert}
\end{equation}

$\kappa_{f}(A,x)$ and $\kappa_{b}(A,b)$ are fixed for a given problem, $Ax=b$. Condition number of the matrix is given by the product of $\kappa_{f}(A,x)$ and $\kappa_{b}(A,b)$,

\begin{equation}
\nonumber
    {\kappa} = \lVert A\rVert  \lVert A^{-1}\rVert 
\end{equation}

\subsubsection{Evaluation metrics}

The following is the \textit{relative} error in our estimates of the relative error at the $k^{th}$ iteration:

\begin{equation}\label{eq:1}
 \max \left( \dfrac {  \dfrac{ \chi_{k} }{\lVert x \rVert} - \dfrac{\lVert \epsilon_{k} \rVert}{\lVert x \rVert} }  { \dfrac{\lVert \epsilon_{k} \rVert}{\lVert x \rVert} }, \dfrac { \dfrac{\lVert \epsilon_{k} \rVert}{\lVert x \rVert} - \dfrac{ \chi_{k} }{\lVert x \rVert} } { \dfrac{ \chi_{k} }{\lVert x \rVert} }  \right)
\end{equation}
 
where $\chi_k$ is the estimate of norm of error at $k^{th}$ iteration and $x$ is the solution to the problem $Ax=b$. The maximum of the two terms inside the brackets encapsulate the fact that an error estimate can either be greater than or less than the true relative error at any iteration, and it can effect a corresponding degree of under-computation or an over-computation. Equation \eqref{eq:1} can also be concisely rewritten as, 

\begin{equation} \label{eq:1_brief}
   \dfrac {  \left \lvert  \dfrac{ \chi_{k} }{\lVert x \rVert} - \dfrac{\lVert \epsilon_{k} \rVert}{\lVert x \rVert} \right \rvert } { \min \left( \dfrac{\lVert \epsilon_{k} \rVert}{\lVert x \rVert} ,  \dfrac{ \chi_{k} }{\lVert x \rVert}  \right) } 
\end{equation}

While \eqref{eq:1_brief} quantifies the uncertainty using the relative residue or an error estimator at a particular iteration, \textit{Uncertainty Ratios} and \textit{Linear Uncertainty Ratios} i.e. $U.R.(\chi)$ and $L.U.R.(\chi)$, are defined to represent average uncertainty in estimating true relative error over $n-d$ iterations, as follows:

\begin{equation}\label{UR2_def_est}
    U.R.(\chi) = \dfrac{1}{n-d} \sum_{k = 0}^{n-d-1} \left(    \dfrac {  \left \lvert  \dfrac{ \chi_{k}^2 }{\lVert x \rVert^2} - \dfrac{\lVert \epsilon_{k} \rVert^2}{\lVert x \rVert^2} \right \rvert } { \min \left( \dfrac{\lVert \epsilon_{k} \rVert^2}{\lVert x \rVert^2} ,  \dfrac{ \chi_{k}^2 }{\lVert x \rVert^2}  \right) }  \right )
\end{equation}

\begin{equation}\label{UR1_def_est}
    L.U.R.(\chi) = \dfrac{1}{n-d} \sum_{k = 0}^{n-d-1} \left(    \dfrac {  \left \lvert  \dfrac{ \chi_{k} }{\lVert x \rVert} - \dfrac{\lVert \epsilon_{k} \rVert}{\lVert x \rVert} \right \rvert } { \min \left( \dfrac{\lVert \epsilon_{k} \rVert}{\lVert x \rVert} ,  \dfrac{ \chi_{k} }{\lVert x \rVert}  \right) }  \right )
\end{equation}

where $d$ is the delay in estimation in terms of iterations. When the relative residual is used for stopping the iterative solution of a linear system, given a tolerance in relative error, the uncertainty is correspondingly:

\begin{equation}\label{UR2_def_res}
    U.R.(r) = \dfrac{1}{n-d} \sum_{k = 0}^{n-d-1} \left(    \dfrac {  \left \lvert  \dfrac{ \lVert r_k \rVert^2 }{\lVert b \rVert^2} - \dfrac{\lVert \epsilon_{k} \rVert^2}{\lVert x \rVert^2} \right \rvert } { \min \left( \dfrac{\lVert \epsilon_{k} \rVert^2}{\lVert x \rVert^2} ,  \dfrac{ \lVert r_k \rVert^2 }{\lVert b \rVert^2}  \right) }  \right )
\end{equation}

\begin{equation}\label{UR1_def_res}
    L.U.R.(r) = \dfrac{1}{n-d} \sum_{k = 0}^{n-d-1} \left(    \dfrac {  \left \lvert  \dfrac{ \lVert r_k \rVert }{\lVert b \rVert} - \dfrac{\lVert \epsilon_{k} \rVert}{\lVert x \rVert} \right \rvert } { \min \left( \dfrac{\lVert \epsilon_{k} \rVert}{\lVert x \rVert} ,  \dfrac{ \lVert r_k \rVert }{\lVert b \rVert}  \right) }  \right )
\end{equation}

In our subsequent analysis, we present a theorem on an upperbound of the expectation of $U.R.(r)$ and $L.U.R.(r)$ over converging problems having fixed singular values and constant forward condition number, we show that $\E(U.R.(\chi)) \leq \E(U.R.(r))$, and conclude that use of an error estimator in such cases as a stopping criterion gives increased savings and accuracy in computation on average, owing to their robustness.

\subsubsection{Theorem on Expectation of $ U.R.$ and $L.U.R.$ } 
\begin{theorem}\label{Theorem_UR}
For linear problems of non-singular real matrices A with a set of fixed singular values, fixed forward condition number $\kappa_{f}(A,x) >> 1$, fixed backward condition number $\kappa_{b}(A,b) >> 1$,
\begin{equation}
\nonumber
\dfrac{\alpha}{\kappa_{b}(A,b)^2} +  \dfrac{1-\alpha}{\kappa_{f}(A,x)^2} \leq \E \left( U.R.(r) \right) + o(2) \leq  \alpha\kappa_{f}(A,x)^2 + (1-\alpha) \kappa_{b}(A,b)^2          
\end{equation}
\begin{center}
\&    
\end{center}

\begin{equation}
\nonumber
\dfrac{\alpha}{\kappa_{b}(A,b)} +  \dfrac{1-\alpha}{\kappa_{f}(A,x)} \leq \E \left( L.U.R.(r) \right) + o(2) \leq  \alpha\kappa_{f}(A,x) + (1-\alpha) \kappa_{b}(A,b)          
\end{equation}

where $n$ is the dimension of the matrix $A$ and $\alpha$ is the average probability of relative residual over predicting the relative error of solution i.e. $\alpha = \dfrac{1}{n-d}\sum_{k=0}^{n-d-1} P_{cond}\left( \dfrac{ \lVert r_k \rVert }{\lVert b \rVert} \geq \dfrac{ \lVert \epsilon_k \rVert }{\lVert x \rVert} \right)$, where $P_{cond}$ refers to conditional probability with conditions stated at the start of the theorem. The sample space for conditional probability and conditional expectation is generated by sampling left-singular matrix $U$ and right-singular matrix $V$ from uniformly distributed orthogonal matrices and unit vector $\hat{b}$ from Gaussian distribution. 
\end{theorem}
\begin{proof}

Let the uncertainty ratio at $k^{th}$ iteration be
\begin{equation*}
    U.R.^{(k)}(r) =    \dfrac {  \left \lvert  \dfrac{ \lVert r_k \rVert^2 }{\lVert b \rVert^2} - \dfrac{\lVert \epsilon_{k} \rVert^2}{\lVert x \rVert^2} \right \rvert } { \min \left( \dfrac{\lVert \epsilon_{k} \rVert^2}{\lVert x \rVert^2} ,  \dfrac{ \lVert r_k \rVert^2 }{\lVert b \rVert^2}  \right) }   
\end{equation*}
and 
\begin{equation*}
    \alpha_k = P_{cond}\left( \dfrac{ \lVert r_k \rVert }{\lVert b \rVert} \geq \dfrac{ \lVert \epsilon_k \rVert }{\lVert x \rVert} \right)
\end{equation*}

Using the definition in \eqref{eq:1} for $  U.R.^{(k)}(r) $, we can write,
\[
    \E \left(  U.R.^{(k)}(r)  \right) = \alpha_k \E_1 \left( \dfrac { \dfrac{ \lVert r_k \rVert^2 }{\lVert b \rVert^2} - \dfrac{\lVert \epsilon_{k} \rVert^2}{\lVert x \rVert^2} } { \dfrac{\lVert \epsilon_{k} \rVert^2}{\lVert x \rVert^2} }   \right) + (1-\alpha_k) \E_2\left( \dfrac { \dfrac{\lVert \epsilon_{k} \rVert^2}{\lVert x \rVert^2} - \dfrac{ \lVert r_k \rVert^2 }{\lVert b \rVert^2}  } { \dfrac{ \lVert r_k \rVert^2 }{\lVert b \rVert^2} }   \right) 
\]

where $\E_1$ refers to $\E\left(  \bullet  \Biggr | \dfrac{ \lVert r_k \rVert }{\lVert b \rVert} \geq \dfrac{ \lVert \epsilon_k \rVert }{\lVert x \rVert} \right)$ and $\E_2$ refers to $\E\left(  \bullet  \Biggr | \dfrac{ \lVert r_k \rVert }{\lVert b \rVert} \leq \dfrac{ \lVert \epsilon_k \rVert }{\lVert x \rVert} \right)$.

\[
     \E \left(  U.R.^{(k)}(r)  \right) = \alpha_k \E_1 \left(  \dfrac{ \lVert r_k \rVert^2 }{\lVert b \rVert^2} \dfrac{\lVert x \rVert^2}{\lVert \epsilon_{k} \rVert^2} -1  \right) + \left(1-\alpha_k\right) \E_2 \left(   \dfrac{\lVert \epsilon_{k} \rVert^2}{\lVert x \rVert^2} \dfrac{ \lVert  b \rVert^2 }{\lVert r_k \rVert^2} -1  \right)
\]

\[
    \E \left(  U.R.^{(k)}(r)  \right) = \alpha_k  \E_1  \left(  \dfrac{ \lVert r_k \rVert^2 }{\lVert \epsilon_{k} \rVert^2}  \dfrac{\lVert x \rVert^2}{\lVert b \rVert^2}  \right) + \left(1-\alpha_k\right)  \E_2 \left(   \dfrac{\lVert \epsilon_{k} \rVert^2}{\lVert r_k \rVert^2} \dfrac{\lVert b \rVert^2}{\lVert x \rVert^2}  \right) -1
\]

Consider $S_k= \dfrac{ \lVert r_k \rVert }{\lVert \epsilon_{k} \rVert}  \dfrac{\lVert x \rVert}{\lVert b \rVert} $, thus, $S_k \in \left \{ \dfrac{1}{\kappa_{b}(A,b)} , \kappa_{f}(A,x) \right \}$,

\[
    \E \left(  U.R.^{(k)}(r)  \right) = \alpha_k  \E  \left(  S_k^2 \big | S_k \geq 1  \right) + \left(1-\alpha_k\right)   \E  \left( \dfrac{1}{S_k^2}  \big | S_k \leq 1  \right) -1
\]

\[
    \E \left(  U.R.^{(k)}(r)  \right) = \alpha_k \left( \E(S_k^2) - o(1) \right) + (1-\alpha_k) \left( \E\left(\dfrac{1}{S_k^2}\right) - o(1) \right) - 1
\]

\[
    \E \left(  U.R.^{(k)}(r)  \right) = \alpha_k \left( \E(S_k^2) \right) + (1-\alpha_k) \left( \E\left(\dfrac{1}{S_k^2}\right) \right) - o(2)
\]

\[
    \E\left( S_k^2 \right) = \E \left( \dfrac{ \lVert r_k \rVert^2 }{\lVert \epsilon_{k} \rVert^2}  \dfrac{\lVert x \rVert^2}{\lVert b \rVert^2} \right)  
\]

Based on fixed conditions of expectation, $\dfrac{\lVert x \rVert}{\lVert b \rVert}$ is constant,

\[
    \E\left( S_k^2 \right) =  \dfrac{\lVert x \rVert^2}{\lVert b \rVert^2} \E \left( \dfrac{ \lVert r_k \rVert^2 }{\lVert \epsilon_{k} \rVert^2}   \right) = \dfrac{\lVert x \rVert^2}{\lVert b \rVert^2} \E \left( \dfrac{ \lVert A\epsilon_{k} \rVert^2 }{\lVert \epsilon_{k} \rVert^2}   \right)    
\]
By Using definition of $l_2$ norm of matrix,
\[
    \dfrac{\kappa_{f}(A,x)^2}{\kappa(A)^2}\leq \E\left( S_k^2 \right) \leq \kappa_{f}(A,x)^2
\]

Similarly,
\[
   \dfrac{\kappa_{b}(A,b)^2}{\kappa(A)^2}\leq  \E\left( \dfrac{1}{S_k^2} \right) \leq \kappa_{b}(A,b)^2
\]

Hence, using above bounds and $\alpha$ as defined in theorem \ref{Theorem_UR},

\[
\dfrac{\alpha}{\kappa_{b}(A,b)^2} +  \dfrac{1-\alpha}{\kappa_{f}(A,x)^2} - o(2) \leq \E \left( U.R.(r) \right) \leq \alpha\kappa_{f}(A,x)^2 + (1-\alpha) \kappa_{b}(A,b)^2 - o(2)     
\]

\begin{equation}
\nonumber
\dfrac{\alpha}{\kappa_{b}(A,b)^2} +  \dfrac{1-\alpha}{\kappa_{f}(A,x)^2} \leq \E \left( U.R.(r) \right) + o(2) \leq  \alpha\kappa_{f}(A,x)^2 + (1-\alpha) \kappa_{b}(A,b)^2          
\end{equation}

Similarly, we can derive the following bounds for $L.U.R(r)$:

\begin{equation}
\nonumber
\dfrac{\alpha}{\kappa_{b}(A,b)} +  \dfrac{1-\alpha}{\kappa_{f}(A,x)} \leq \E \left( L.U.R.(r) \right) + o(2) \leq  \alpha\kappa_{f}(A,x) + (1-\alpha) \kappa_{b}(A,b)    
\end{equation}

\end{proof}

\subsubsection{Note on Delay-based Error Estimators } 

For delay-based estimators, $\chi_k^2 \approx \norm{\epsilon_k}^2 - \norm{\epsilon_{k+d}}^2$, and 
\[
      U.R.^{(k)}(\chi) \approx \dfrac{\norm{\epsilon_{k+d}}^2}{\norm{\epsilon_k}^2 - \norm{\epsilon_{k+d}}^2}.
\]
On the other hand, the relative residual has a strict lower bound of $\dfrac{1}{\kappa(A)} \dfrac{\norm{\epsilon}}{\norm{x}}$ and a strict upper bound of $\kappa(A) \dfrac{\norm{\epsilon}}{\norm{x}}$, %$\left( \dfrac{1}{\kappa(A)} \dfrac{\norm{\epsilon}}{\norm{x}} \leq \dfrac{\norm{r}}{\norm{b}} \leq  \kappa(A) \dfrac{\norm{\epsilon}}{\norm{x}} \right)$ 
for a given matrix $A$. this allowed us to provide conditional upper and lower bounds when $\norm{\epsilon_{k+d}}^2 < \norm{\epsilon_{k}}^2$. Hence, use of delay based error estimators is reliable when an iterative method like GMRES or BiCG has a decreasing $l_2$ norm of error for most iterations. It can be argued that value of $d$ can be increased to an extent such that $\norm{\epsilon_{k+d}}^2 < \norm{\epsilon_{k}}^2$ for all the iterations, but in practice, a constraint ($d \ll n$) is required due to low latency requirements. The above conclusions are valid without any restriction on the dimension of the matrices.

\section{Numerical results}\label{Section4}

\subsection{Modified $l_2$ norm error estimator for GMRES}

In this experiment, we empirically highlight the advantages of including term-2 in the approximation of the GMRES error, as proposed in equation \eqref{newesteq}. We conducted experiments on non-symmetric positive definite high condition number matrices (generated via eigen decomposition with eigen values spread exponentially between $1$ and $\kappa(A)$) in which we compared the values of different error estimators and relative residue on each iteration.

\begin{centering}

\begin{figure}[]
      \includegraphics[width=0.8\linewidth, keepaspectratio]{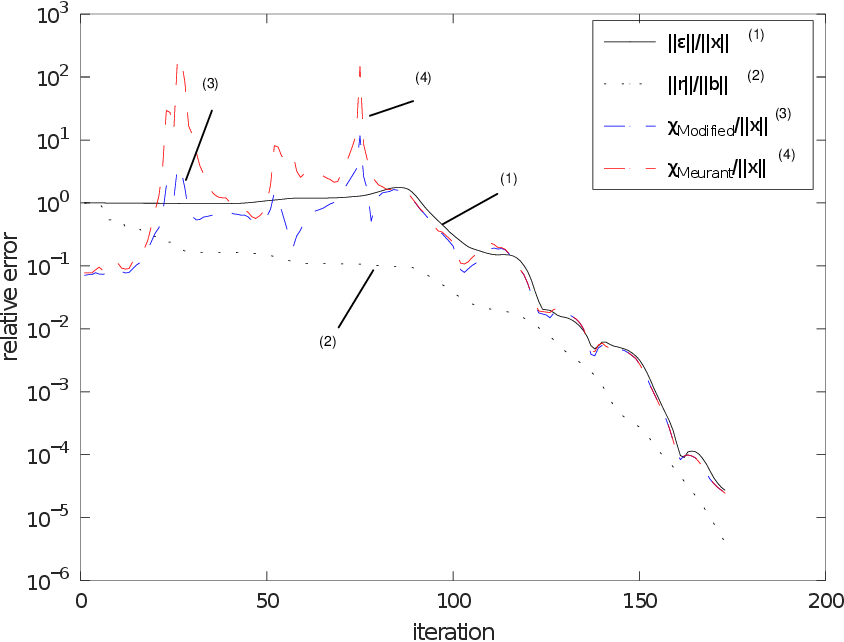}
    \caption{Convergence plot showing comparison of different relative error estimates and true relative error for randomly generated non-symmetric matrix ($\kappa(A) \approx 10^8$) with random right hand side vector $b$, and the delay in iterations $d$ = 10. The matrix has a dimension of 500 and its first 20 largest eigen values have a magnitude $\approx 10^7$ and the rest in the range $[0.1,10]$. Note that the modified estimator is much closer to true error than the original estimator as it avoids unstable overshoots seen in the absence of the offset term, and the logarithmic scale should be noted. $L.U.R$ (residual) = 12.1, $L.U.R$ (Meurant) = 10.1, and $L.U.R$ (Modified) = 1.2. }
    \label{fig:ust2}
\end{figure}

\end{centering}

In Figure \ref{fig:ust2}, we compare relative $l_2$ norm of error estimator proposed by Meurant \cite{meurant2011estimates}, the modified $l_2$ norm error estimator proposed in equation \eqref{newesteq}, the relative residual and the true relative error at each iteration. We observe that including term-2 in equation \eqref{newesteq} helped in avoiding unstable overshoots in estimation for some of the iterations as seen in Figure \ref{fig:ust2}.

\begin{figure}[]
      \includegraphics[width=0.8\linewidth, keepaspectratio]{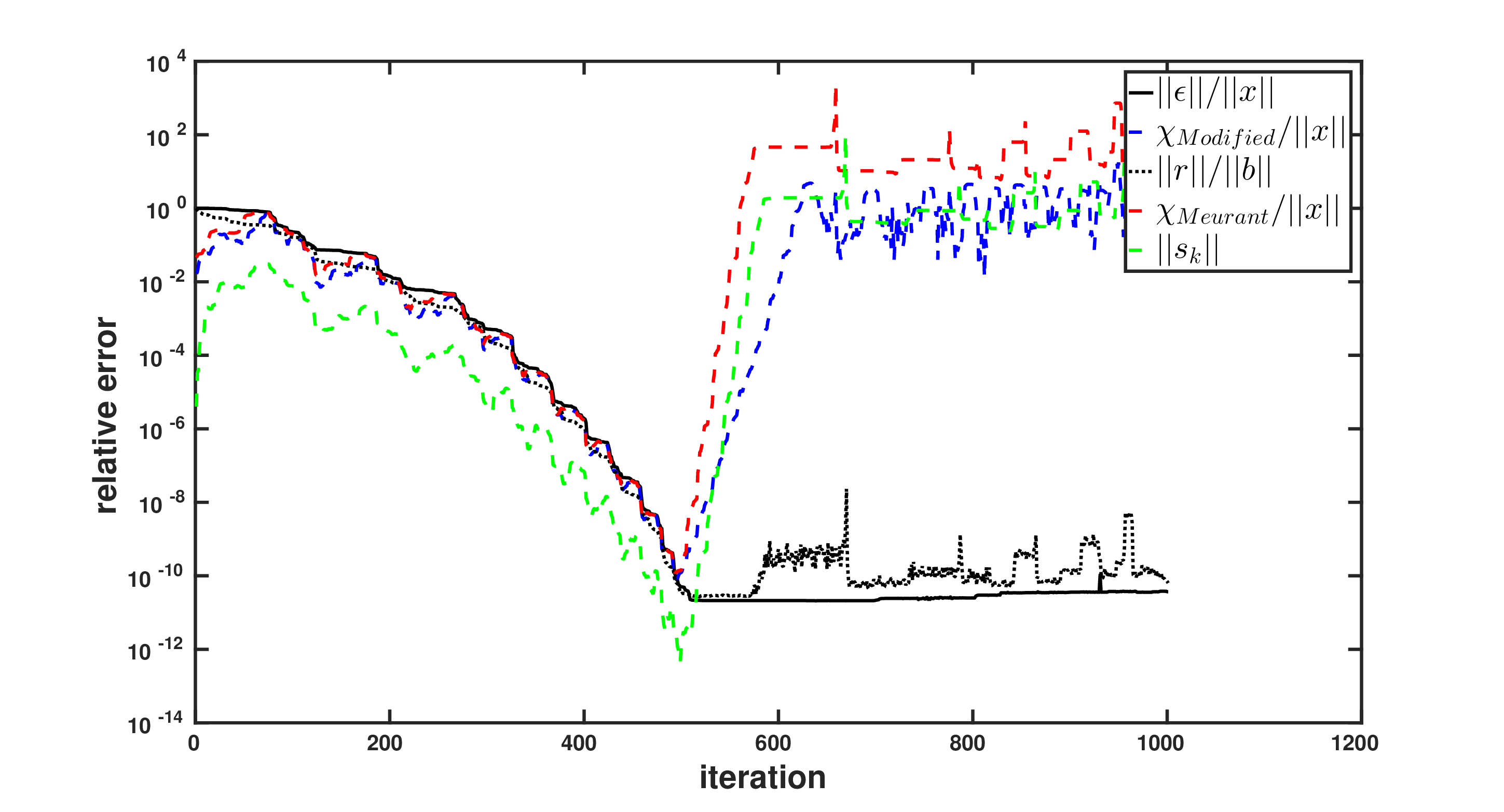}
    \caption{Behaviour of relative error estimates in GMRES (for a symmetric positive definite matrix with $\kappa(A) = 10^6$) when the numerical precision is exhausted. Note that $\norm{s_k}$ and the two error estimators blow up while the actual error and residue are indeed small. The term $\norm{s_k}$ can also be used as a trigger for detecting the exhaustion of numerical precision.}
    \label{blowup}
\end{figure}

It is also observed that GMRES error estimators may behave erratically when numerical precision of computing system is exhausted as seen in Figure \ref{blowup}. The exhaustion of numerical precision can lead to near singularity of Hessenberg matrix $H_k$ formed during the Arnoldi iteration which can cause large errors in $H_k^{-1}$ and its functions. Note that $\norm{s_k}$ is also a function of $H_k^{-1}$ and can be used as trigger to predict this exhaustion of numerical precision as seen in Figure \ref{blowup} making this error estimate a robust stopping criterion.

\subsection{A-norm and $l_2$ norm error estimators for Bi-CG}

In this experiment, we want to empirically evaluate the behaviour of delay-based BiCG error estimators derived in section \ref{BiCGSection} when compared to relative residue, in problems with a high condition number. We considered non-symmetric matrices from the Matrix Market\footnote{https://math.nist.gov/MatrixMarket/} dataset ('sherman5') with a random right hand side vector $b$. 

\begin{figure}[]
  \centering
  \includegraphics[width=0.8\linewidth]{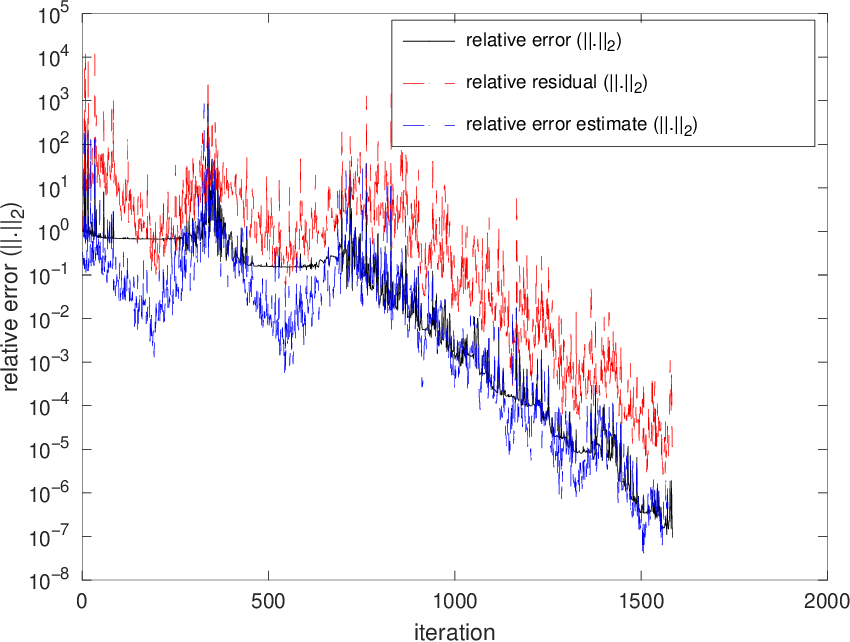}
  \caption{Estimation of relative $l_2$-norm of error for a Non-symmetric matrix (`sherman5') from the Matrix Market dataset with a random right hand side vector, using the BiCG error estimator and the relative residue; dimension of the matrix = $3312 \times 3312$; condition number of the matrix is $3.9 \times 10^{5}$ and $d = 10$.}
  \label{fig:figure 3}
\end{figure}

\begin{figure}[H]
  \centering
  \includegraphics[width=0.8\linewidth]{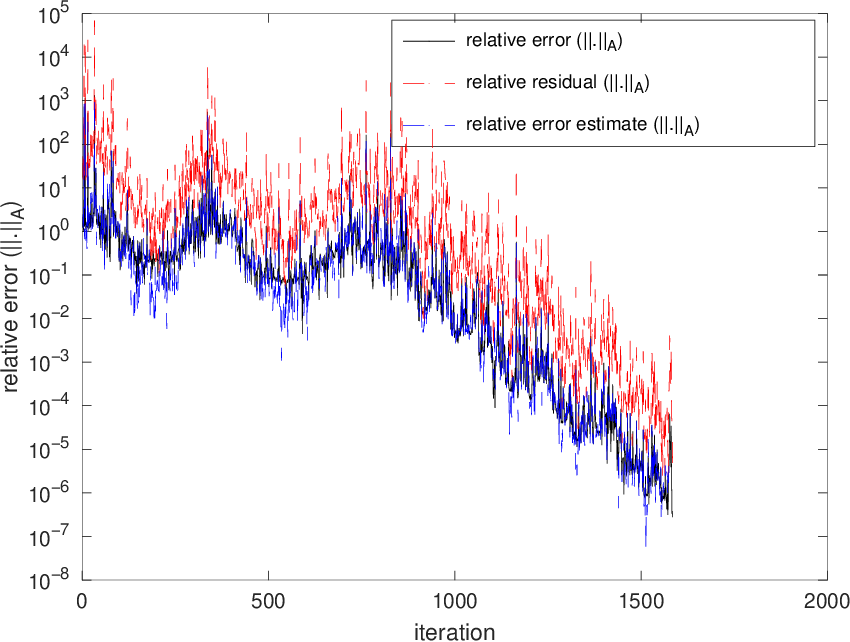}
  \caption{Estimation of relative A-norm of error for a Non-symmetric square matrix (`sherman5') from the Matrix Market dataset with a random right hand side vector, using the BiCG error estimator and the relative residue; absolute values are considered for $r^{T}A^{-1}r$ and its approximation; dimension of the matrix = $3312 \times 3312$; condition number of the matrix is $3.9 \times 10^{5}$ and $d = 10$.}
  \label{fig:figure 2}
\end{figure}

In Figure \ref{fig:figure 2}, we compare estimated relative A-norm of the error, the relative residual, and the true relative A-norm of the error at each iteration. Figure \ref{fig:figure 3} shows the comparison between estimated relative $l_{2}$ norm error, the true relative $l_{2}$ norm of the error, and the relative residual. We observe that the proposed BiCG error estimators performs consistently across iterations in both cases, where the error estimator is reliable as an indicator of the actual error for stopping. 

\subsection{Experimental results on performance of error estimators on randomized set of matrices} \label{expt_errorestimate}

The objective of the experiments\footnote{Code freely available at https://github.com/CSPL-IISc/Error-Linsolve} is to asses the predictive performance of error estimators on large randomized set of matrices. This study is useful in understanding the variability in average behaviour of estimators with respect to different aspects of problem. In previous sections, we had defined the Uncertainty Ratio as the average error in any of the predictors of true relative error that may be used for stopping.

\begin{figure}[]
    \centering
    \subfloat[$L.U.R$ vs. $\kappa_{f}(A,x)$]{\includegraphics[width = 1\textwidth]{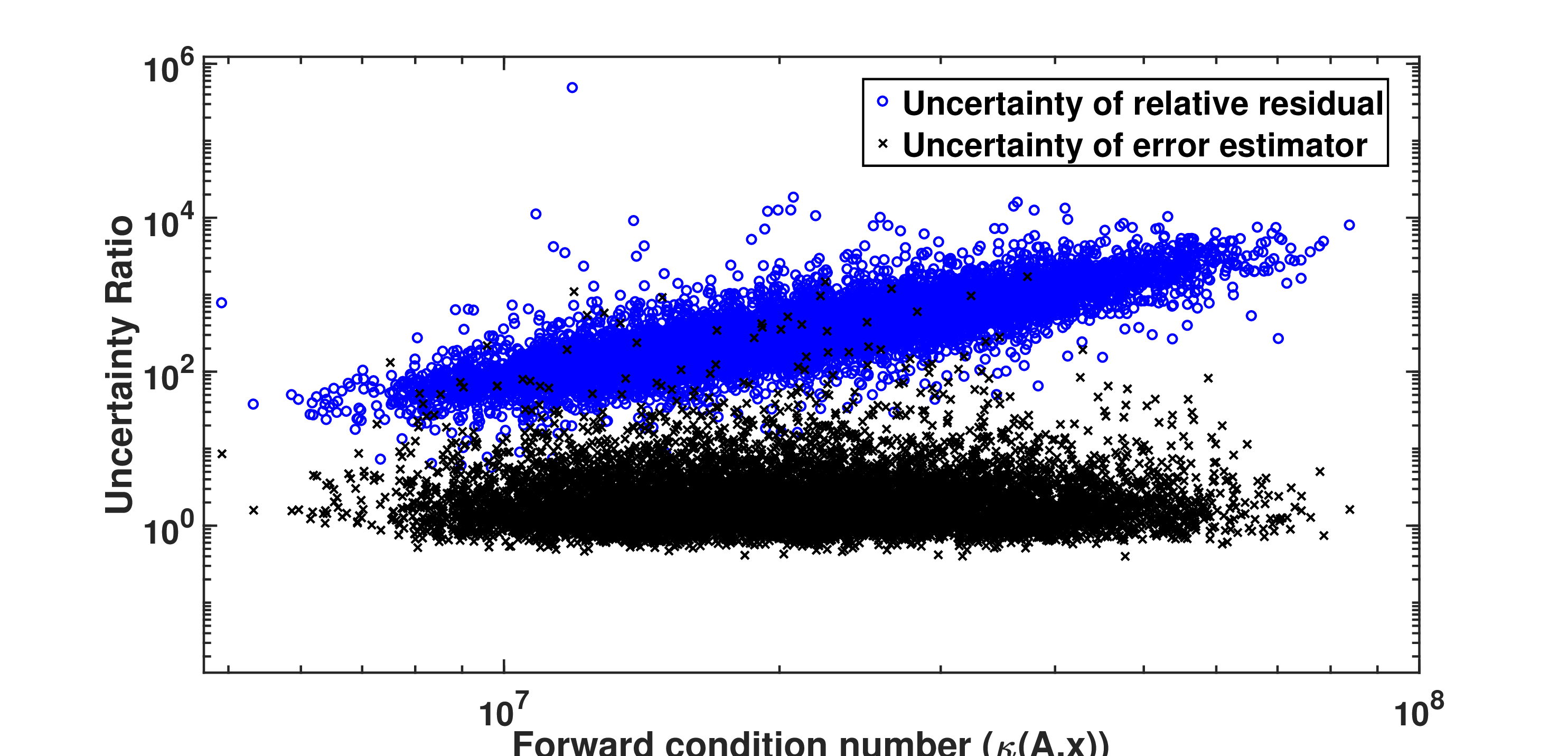}}\\
     \subfloat[$L.U.R$ vs. $\kappa_{b}(A,b)$]{\includegraphics[width = 1\textwidth]{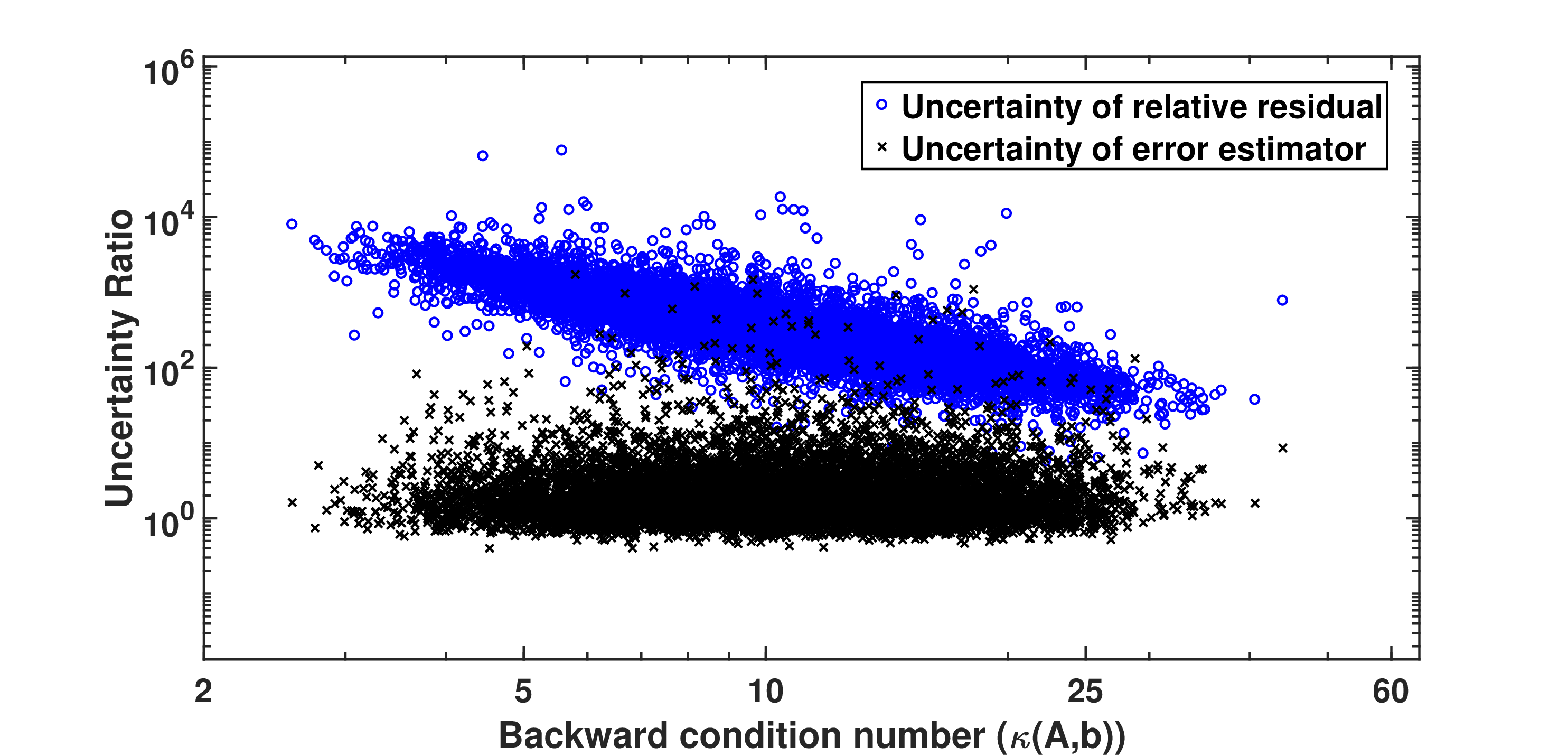}} 
    \caption{Uncertainty $L.U.R$ when using the estimator for the relative $l_2$-norm of error, or the relative residue, for stopping the BiCG algorithm when a tolerance in relative error was given. The forward and backward condition number of problems were varied. Mean $L.U.R$(residual) = 288. Mean $L.U.R$(estimator) = 5.9.}
    \label{UR_100_BiCG}
\end{figure}

In this experiment, we measure the Uncertainty Ratio of estimator and residual for a given forward condition number and backward condition number for each problem in randomized set. The randomized dataset was chosen to consist of 10,000 problems of high condition numbers of dimension $100$ on which the evaluation metrics $L.U.R(r)$ and $L.U.R(\chi)$ are measured with $d=10$ for delay in estimation. We consider two kinds of non-symmetric matrices i.e. positive definite and general matrices in the dataset. The right hand side vectors $b$ are random vectors from a unit normal distribution of same dimension. 

\begin{figure}[]
    \centering
    \subfloat[$L.U.R$ vs. $\kappa_{f}(A,x)$]{\includegraphics[width = 1\textwidth]{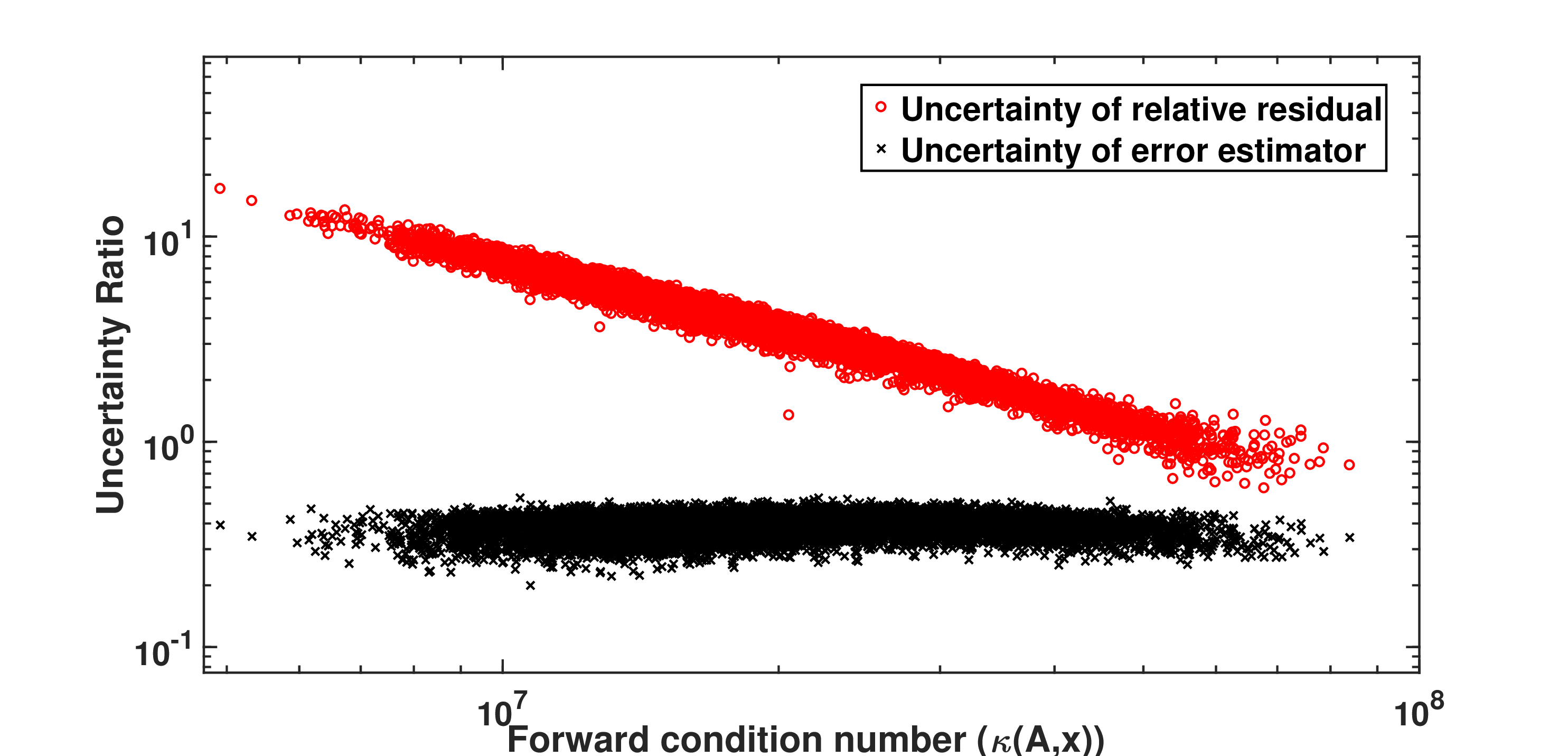}} \\
     \subfloat[$L.U.R$ vs. $\kappa_{b}(A,b)$]{\includegraphics[width = 1\textwidth]{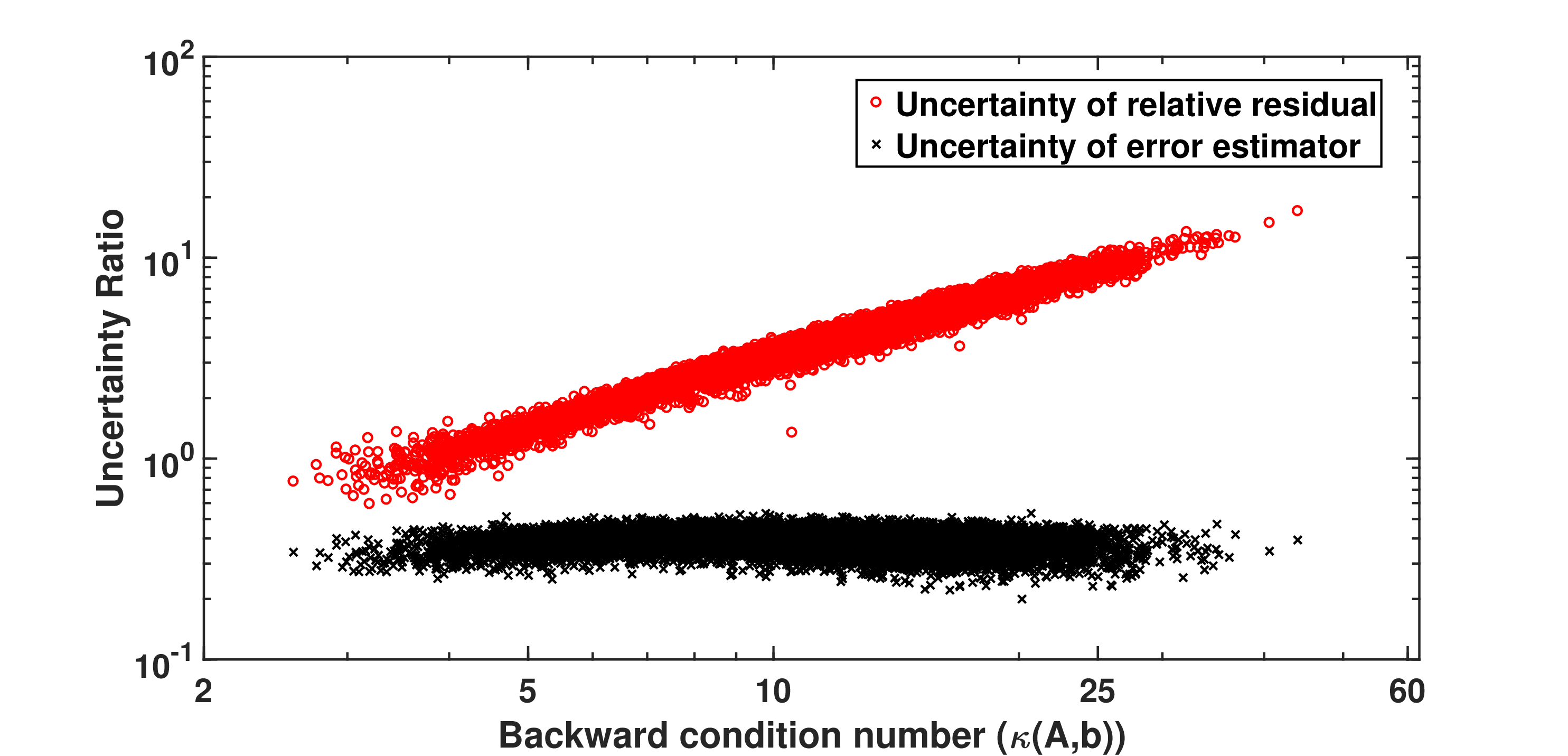}} 
    \caption{Uncertainty $L.U.R$ when using the estimator for the relative $l_2$-norm of error, or the relative residue, for stopping the GMRES algorithm when a tolerance in relative error was given. The forward and backward condition number of problems were varied. Mean $L.U.R$ (residual) = 2.49. Mean $L.U.R$ (estimator) = 0.286. Note that the modified and original GMRES estimators are nearly identical in these scatter plots where we have only a small number of dimensions n=100 in the problems.}
    \label{UR_100_GMRES}
\end{figure}

Figure \ref{UR_100_BiCG} and \ref{UR_100_GMRES} show the behaviour of Linear Uncertainty Ratio ($L.U.R.$), defined by \eqref{UR1_def_est} and \eqref{UR1_def_res}, with the condition number of the problem (both forward and backward) for BiCG and GMRES algorithms respectively. Here each circle (blue for BiCG and red for GMRES) represents the uncertainty ratio corresponding to the relative residual for a particular problem i.e for a particular matrix $A$ and a right hand side vector $b$ with $x_0$ sampled from standard Gaussian distribution. Similarly, each (black) cross represents the uncertainty ratio corresponding to the relative error estimator.

\begin{figure}[]
    \centering
    \subfloat[BiCG ($L.U.R.$ vs. $\kappa_{f}(A,x)$)]{\includegraphics[width = 1\textwidth]{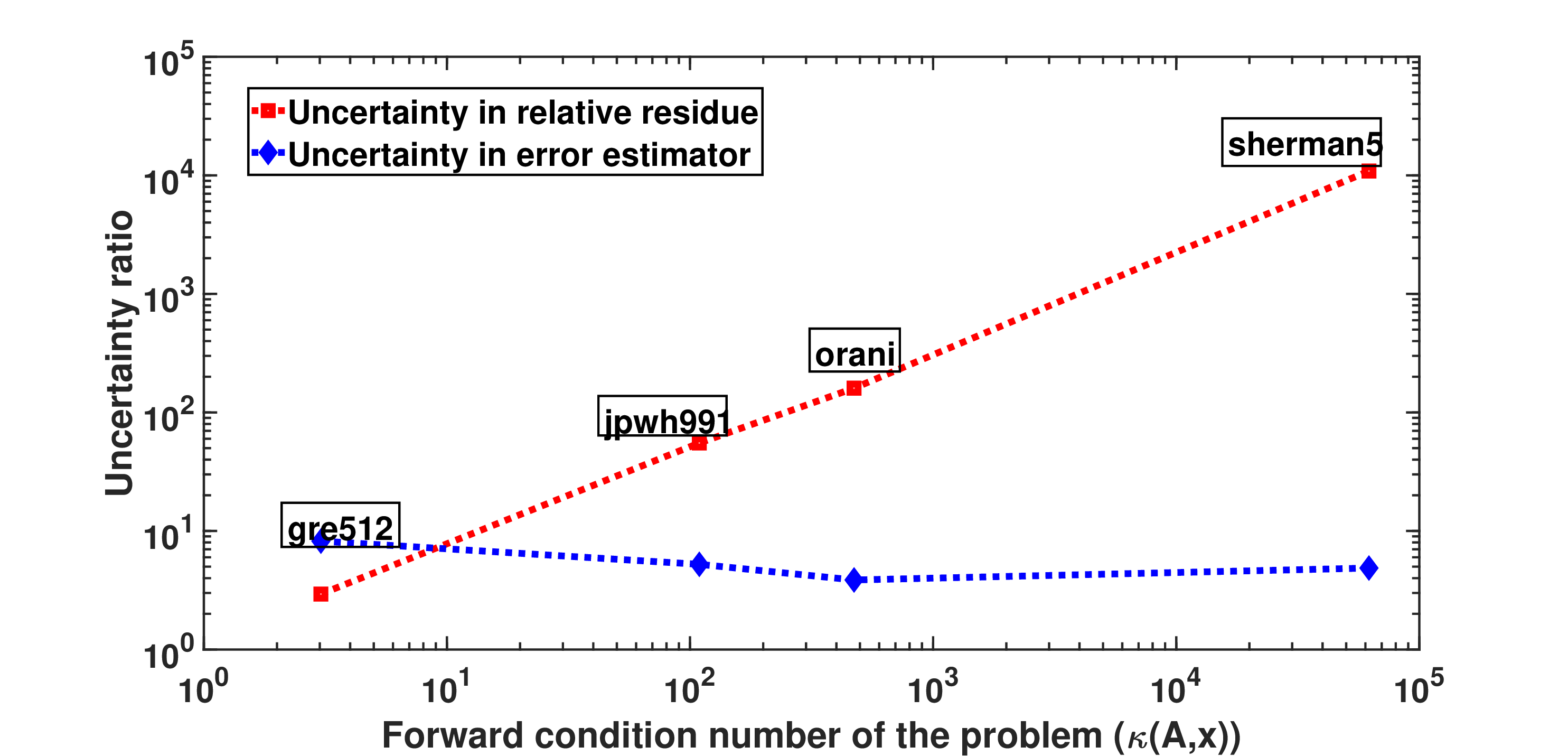}} \\
     \subfloat[GMRES ($L.U.R.$ vs. $\kappa_{b}(A,b)$)]{\includegraphics[width = 1\textwidth]{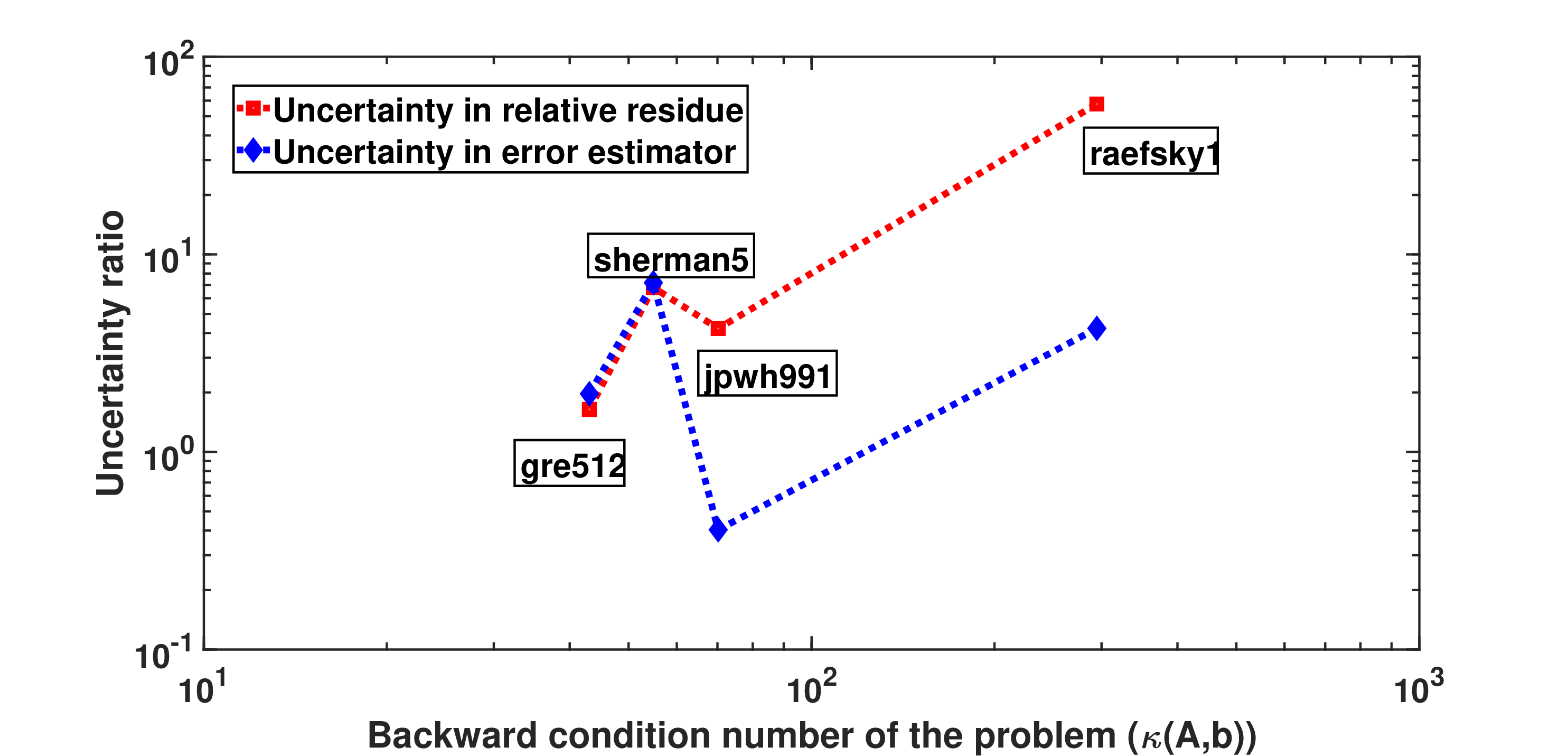}} 
    \caption{Uncertainty $L.U.R$ when using the estimator for the relative $l_2$-norm of error, or the relative residue, for stopping the iterations when a tolerance in relative error was given; matrices (with sizes) considered from the Matrix-Market data set for this evaluation were: `gre15' (512), `jpwh991'(991), `orani' (2529), `sherman5'(3312), and `raefsky1' (3242).}
    \label{UR_matrixmarket}
\end{figure}

We observe that there is a large increase in L.U.R.(r) of relative residual in high condition number problems, whereas the L.U.R.($\chi$) for the relative error estimator remains bounded, indicating the robustness of BiCG and GMRES error estimators. It is also interesting to note that Uncertainty Ratio is dependent on both forward and backward condition number of the problem. For BiCG, the $L.U.R.(r)$ increases with the forward condition number of the problem, while $L.U.R.(r)$ increases with the backward condition number in case of the GMRES algorithm. Note that the $\alpha$ values in the result of the theorem, that determine the increasing (or decreasing) relationship with the forward (or backward) condition number, can vary across algorithms. Since the same dataset was used for testing the metrics on both Bi-CG and GMRES algorithms, the scales on x-axis differ in the two figures (and note $\kappa_{f}(A,x) * \kappa_{b}(A,b) = \kappa$). While Bi-CG has much larger uncertainty ratios because of its oscillatory covergence behavior, the
expected under and over computation for a required tolerance in error are comparable for both Bi-CG and GMRES, as shown in the next section.

For further demonstration of the above conclusions even in matrices of practical applications, non-symmetric matrices of varied condition numbers and dimensions were collected from Matrix Market and Matrix Sparse Collection data set and were tested with corresponding random normalized right hand side vectors. Figure \ref{UR_matrixmarket} shows the relation between Uncertainty Ratio and condition numbers (forward and backward) of the corresponding problems. We observe the significant increase in the uncertainty ratio of the relative residual with the increase in the condition number of the problem while the uncertainty ratio of the error estimator (for both BiCG and GMRES) stays almost stagnant as the condition number of the problem increases, depicting error estimators as a robust stopping criteria. Though we mostly focus on the evaluation metrics for BiCG and GMRES algorithms, one should keep into account that the performance of the CG estimator is implicitly described by the Bi-CG estimator in the case of symmetric positive definite matrices.

\subsection{Under and over computation with BiCG and GMRES estimators} 

While using iterative algorithms, one would like to stop the computation at error tolerances depending on the requirements of the applications. Since there is a larger uncertainty when stopping using the relative residue, as shown in Figures \ref{UR_100_BiCG}, \ref{UR_100_GMRES} and \ref{UR_matrixmarket}, an under or over computation is expected for any desired error tolerance in solution. One might stop the iterations too early (under computation) or deploy too many iterations (over computation) due to this uncertainty in error of the solution, resulting in loss of efficiency or accuracy. In this experiment, the objective is to evaluate the degree of over or under computation due to the uncertainty in error when using the relative residual as an indicator for stopping.

\begin{figure}
    \centering
    \subfloat[BiCG]{\includegraphics[width = 1\textwidth]{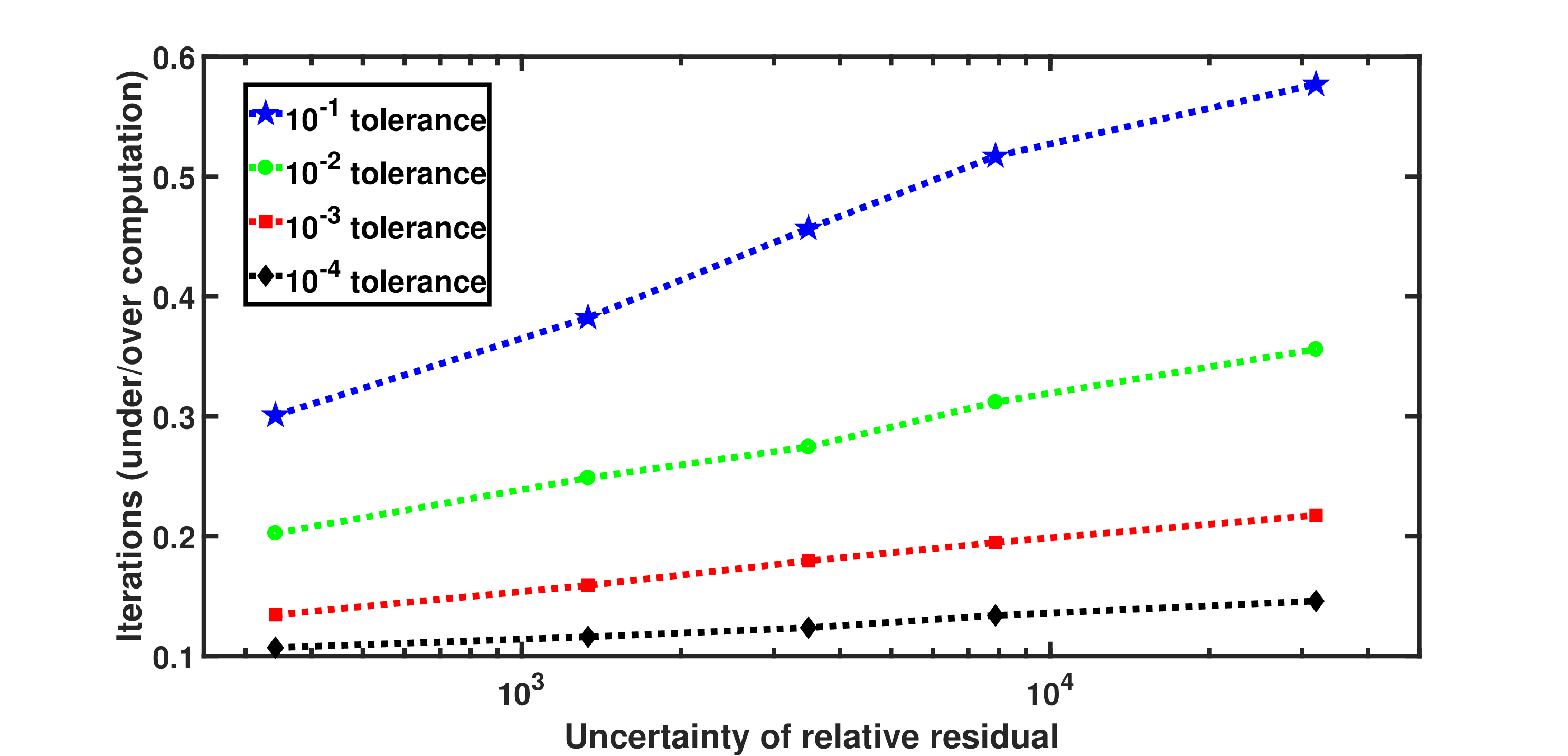}}\\ 
     \subfloat[GMRES]{\includegraphics[width = 1\textwidth]{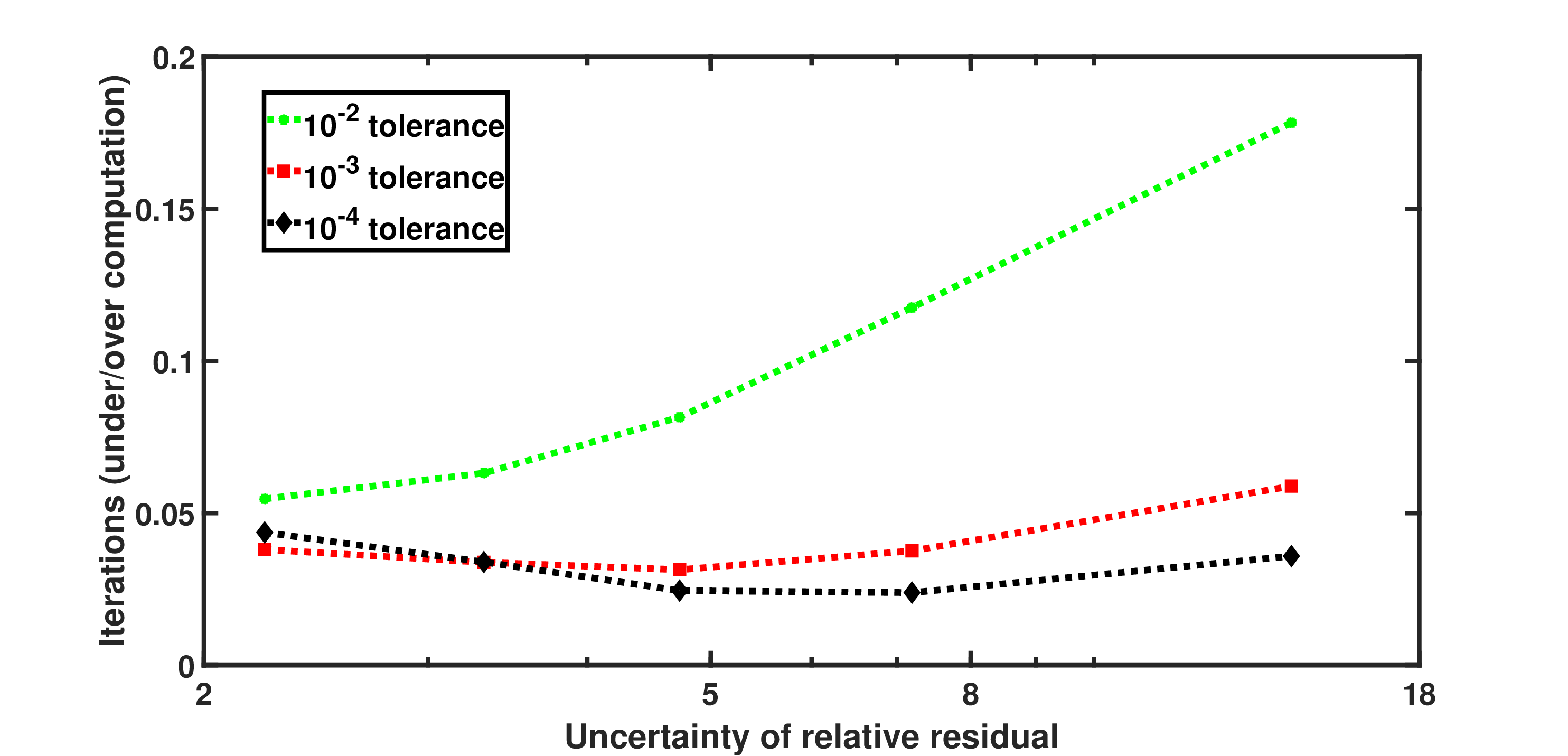}} 
    \caption{Fraction of iterations under (or over) computed when iterations are stopped at different tolerances, as the uncertainty of relative residue ($L.U.R.$) increases. Note that the uncertainty in Bi-CG can be larger due to oscillatory convergence, but the expected over and under computation in using the residue for stopping, is comparable for Bi-CG and GMRES.  The forward and backward condition numbers were increased correspondingly for a higher uncertainty.}
    \label{Iterations_100}
\end{figure}

We define the relative iterations as $I$ which are either saved or lost for a desired error tolerance as:

\begin{equation}
\nonumber
    I = \frac{\lvert Iter_{\frac{\norm{r}}{\norm{b}}} - Iter_{\frac{\norm{\chi}}{\norm{x}}} \rvert}{ \min \left( Iter_{\frac{\norm{r}}{\norm{b}}} ,  Iter_{\frac{\norm{\chi}}{\norm{x}}}  \right) }
\end{equation}\label{iterations_saved_lost}

For any tolerance such as $10^{-4}$, $I$ signifies the fraction of iterations over or under computed while using relative residual as a stopping criteria as compared to while using error estimator as stopping criteria. $Iter_{\frac{\norm{r}}{\norm{b}}}$ is the number of iterations required by the relative residual ($\frac{\norm{r}}{\norm{b}}$) to reach the desired tolerance, and we have a similar variable for the iterations required by the error estimate. Figure \ref{Iterations_100} shows the plot of $I$ at error tolerances of $10^{-1},10^{-2},10^{-3}$ and $10^{-4}$ with increasing $L.U.R.(r)$ averaged over the matrix problems used in Section \ref{expt_errorestimate} for both the algorithms. It should be noted that the value of $I = 1$ means there is a 100 percent difference in the iterations required by $\frac{\norm{r}}{\norm{b}}$
and $\frac{\norm{\chi}}{\norm{x}}$ to reach a desired tolerance. We observe the increase in $I$ with the uncertainty ratio of relative residual which signifies the importance of error estimators while using iterative algorithms.

\section{Conclusions}

The significance of error estimators for efficient stopping (or restarting) are clearly evident for problems with condition numbers $\kappa \gtrsim 100$, and is emphasized by numerical examples and an analysis of the expected uncertainty in convergence when using the norm of the residual vector for stopping, for any given tolerance in error. Results showing the reduction of uncertainty ratio while using error estimators were presented for both Bi-CG and GMRES algorithms. It was also highlighted that the proposed estimator for Bi-CG converges to the earlier proposed estimator for the CG algorithm in the case of symmetric positive-definite matrices. Further, a simple analysis of the change in computation due to the estimators was presented along with numerical results of large number of problems of manageable dimensions. While Bi-CG has much larger uncertainty ratios because of its oscillatory behavior (and this can be improved by using extended versions of the algorithm \cite{fletcher1976conjugate} \cite{sleijpen1994bicgstab} \cite{sleijpen1993bicgstab} \cite{van1992bi} \cite{freund1991qmr}), the expected under and over computation for a required tolerance in error are comparable for both Bi-CG and GMRES, as expected. Based on the results discussed in the previous sections, we conclude that the estimate for the A-norm or the $l_{2}$ norm of the error should be implemented into software realization of iterative solvers, instead of using only the relative norm of the residual vector as a criterion for stopping.

\section*{Author declarations:}
No funding was received for conducting this study. The authors have no relevant financial or non-financial interests to disclose.

%\section*{Appendix}
\appendix
\section{Bi-CG and its relation to non-symmetric Lanczos algorithm}\label{BICG}

Bi-Conjugate Gradient (Bi-CG) algorithm is an extension to the Conjugate Gradient (CG) algorithm which is used to solve a system of linear equations and it is suited for a non-symmetric matrix. 

\begin{algorithm}

\caption{Bi-Conjugate Gradient Algorithm }

\label{BiCG algorithm}

\begin{algorithmic}[1]

\Procedure{}{} 
 
  \textbf{input} A, $A^{T}$, b, $x_{0}$,  $y_{0}$ 
  \State$r_{0}$ = b - A$x_{0}$ 
  \State$\tilde{r_{0}}$ = b - $A^{T}$ $y_{0}$ 
  \State$p_{0}$ = $r_{0}$ 
  \State$q_{0}$ = $\tilde{r_{0}}$ 
  \For{k = 1...until convergence} 
 \State ${\alpha}_{k-1}$ = ${\frac {\tilde{r}_{k-1}^T r_{k-1}}{q_{k-1}^T A p_{k-1}}}$ 
  \State$x_{k}$ = $x_{k-1}$ + $\alpha_{k-1}p_{k-1}$ 
  \State$y_{k}$ = $y_{k-1}$ + $\alpha_{k-1}q_{k-1}$ 
  \State $ r_k = r_{k-1} - \alpha_{k-1}Ap_{k-1}$
   \State $ \tilde{r}_k = \tilde{r}_{k-1} - \alpha_{k-1}A^{T}q_{k-1}$
  \State${\beta}_{k-1}$ = ${\frac {\tilde{r_{k}}^T r_{k}}{\tilde{r}_{k-1}^T r_{k-1}}}$ 
  \State$p_{k} = r_{k} + \beta_{k-1}p_{k-1}$ 
 \State $q_{k} = \tilde{r_{k}} + \beta_{k-1}q_{k-1}$ 

\EndFor
\EndProcedure

\end{algorithmic}

\end{algorithm}

Bi-CG can also be derived from the non-symmetric Lanczos algorithm. Let $v_{1}$ and $\tilde{v_{1}}$ be the given starting vectors to the non-symmetric Lancozs algorithm (such that $\lVert v_{1} \rVert = 1$ and $(v_{1},\tilde{v_{1}}) = 1)$, the two three term recurrences which help in forming two bi-orthogonal subspaces can be as follows:

For k=1,2 $\dots$

\begin{equation}
\nonumber
\begin{array}{ccl}
  z_{k} = Av_{k} - w_{k}v_{k} - \eta_{k-1}v_{k-1}\\
  \tilde{z_{k}} = A^{T}\tilde{v_{k}} - w_{k}\tilde{v_{k}} - \tilde{\eta}_{k-1}\tilde{v}_{k-1}
  \end{array}
\end{equation}

The coefficient $w_{k}$ being computed as $w_{k} = (\tilde{v_{k}},Av_{k})$. The other coefficients $\eta_{k}$ and $\tilde{\eta_{k}}$ are chosen (provided ($\tilde{z_{k}},v_{k}) = 0$) such that $\eta_{k}\tilde{\eta_{k}} = (\tilde{z_{k}},z_{k})$ and the new vectors at step $k+1$ are given by:

\begin{equation}
\nonumber
\begin{array}{ccl}
v_{k+1} = \dfrac{z_{k}}{{\eta_{k}}}\\ \\
\tilde{v}_{k+1} = \dfrac{\tilde{z_{k}}}{\tilde{\eta}_{k}}
\end{array}
\end{equation}\hfill\break

The relationship between the bi-orthogonal subspaces and matrix $A$ can be written in the form of a non-symmetric tri-diagonal matrix form (under the condition $\tilde{V_k}^T A V_k = T_k$) as:

\begin{equation}
\nonumber
 T_{k} = \left[
    \begin{array}{ccccc}
    \omega_{1} & \eta_{1} & 0 & \cdots& 0\\
    \tilde{\eta_{1}} & \omega_{2} & \eta_{2} & & \vdots\\
    0 & \ddots&\ddots&\ddots& 0\\
    \vdots& & \tilde{\eta}_{k-2} & \omega_{k-1} & \eta_{k-1}\\
    0 & \cdots& 0 & \tilde{\eta}_{k-1} & \omega_{k} 
    \end{array}
    \right] 
\hfill \break
\end{equation}

\section{ Relating BiCG error estimator to CGQL  } \label{BCGtoCG}

The difference between two consecutive A-norm of error (and A-measure for matrices that are not SPD) in case of a Conjugate Gradient algorithm at iteration $k$ and $k+1$ can be given by: 

\begin{equation}
\nonumber
    \lVert \epsilon_{k} \rVert_{A}^{2}  - \lVert \epsilon_{k+1} \rVert_{A}^{2}  =  \alpha_{k} r_{k}^{T} r_{k} 
\end{equation}

Hence by inducing a delay of $d$ iterations we can easily compute A-norm of error.

 \begin{equation}
 \nonumber
 \begin{array}{ccl}
    \lVert \epsilon_{k} \rVert_{A}^{2} - \lVert \epsilon_{k+d} \rVert_{A}^{2}  =  \sum_{j=k}^{k+d} \alpha_{j} r_{j}^{T} r_{j} \\
    \lVert \epsilon_{k} \rVert_{A}^{2} \approx \sum_{j=k}^{k+d} \alpha_{j} r_{j}^{T} r_{j}
    \end{array}
 \end{equation}

For A-norm estimation of error in BiCGQL algorithm we derived the following results:

\begin{equation}
 \nonumber
    \lVert \epsilon_{k+1} \rVert_{A}^{2} =   r_{k+1}^T A^{-1} r_{k+1} = - \alpha_{k} r_{k}^T p_{k} +   r_{k+1}^T A^{-1} r_{k}   +  \alpha_{k}^{2} p_{k}^T A p_{k}
\end{equation}
  
By using $r_{k+1} = r_k - \alpha_k A p_k$ the above result can also be written as: 
  
\begin{equation}
 \nonumber
   r_{k}^T A^{-1} r_{k} - r_{k+1}^T A^{-1} r_{k+1} = \alpha_{k} r_{k}^T p_{k} +   \alpha_{k} r_{k}^T (A^{T})^{-1} A p_{k}   -  \alpha_{k}^{2} p_{k}^T A p_{k} 
\end{equation} 
   
   For a symmetric matrix $A=A^{T}$ the above equation can be further written as:
   
   \begin{equation}
    \nonumber
        r_{k}^T A^{-1} r_{k} - r_{k+1}^T A^{-1} r_{k+1} = \alpha_{k} r_{k}^T p_{k} +   \alpha_{k} r_{k}^Tp_{k}   -  \alpha_{k}^{2} p_{k}^T A p_{k}
   \end{equation} 
   
Further, using $p_{i}^{T}r_{j} = 0$ for $i \neq j $ in CG, and also $ r_{k+1}  = r_{k} - \alpha_{k} A p_{k} $ and thus substituting $A p_{k} = \dfrac{r_{k} -  r_{k+1}}{\alpha_{k}}$, we reduce the above as:
   
\begin{equation}\label{CGQL_to_BICGQL}
    \begin{array}{ccl}
        r_{k}^T A^{-1} r_{k} - r_{k+1}^T A^{-1} r_{k+1} = \alpha_{k} r_{k}^T p_{k} +   \alpha_{k} r_{k}^Tp_{k}   -  \alpha_{k}p_{k}^Tr_{k} + \alpha_{k}p_{k}^{T}r_{k+1}\\
        \implies r_{k}^T A^{-1} r_{k} - r_{k+1}^T A^{-1} r_{k+1} = \alpha_{k} r_{k}^T p_{k}\\
        \implies r_{k}^T A^{-1} r_{k} - r_{k+1}^T A^{-1} r_{k+1} = \alpha_{k} r_{k}^T r_{k} + \beta_{k-1}r_{k}^{T}p_{k-1}\\
        \implies r_{k}^T A^{-1} r_{k} - r_{k+1}^T A^{-1} r_{k+1} = \alpha_{k} r_{k}^T r_{k} 
    \end{array}
\end{equation} 
   
Equation \ref{CGQL_to_BICGQL} is equivalent to the $A$-norm estimation in CGQL Algorithm \cite{golub1997matrices}, and thus the proposed $A$-norm estimator for BiCG is equivalent to the CGQL estimator when $A=A^T$.

\section{Relations between non-symmetric Lanczos tridiagonal Matrix ($ T_{k} $) and residual vectors ($r$  and $\tilde{r}$) of BiCG algorithm }

A direct relationship between $T_{k}$ , $r_{k}$ and $\tilde{r_{k}}$ can be given by:

  \begin{equation}
       (T_{n}^{-1})_{(1,1)} =   (T_{k}^{-1})_{(1,1)} + \frac { \tilde{r_{k}}^{T}A^{-1}r_{k} } {\lVert r_{0} \rVert^{2} }
       \label{eq: x1}
   \end{equation} 
   
   \begin{equation}
    (T_{n}^{-1})_{(1,1)} =   (T_{k+1}^{-1})_{(1,1)} + \frac { \tilde{r}_{k+1}^{T}A^{-1}r_{k+1} } { \lVert r_{0} \rVert ^{2} }
      \label{eq: x2}
    \end{equation}
    
    Subtracting \ref{eq: x1} and \ref{eq: x2} we get:
    
    \begin{equation}
     \nonumber
        (T_{k+1}^{-1})_{(1,1)} = (T_{k}^{-1})_{(1,1)} + \frac{  ( \tilde{r_{k}}^{T}A^{-1}r_{k} -  \tilde{r}_{k+1}^{T}A^{-1}r_{k+1} ) } { \lVert r_{0} \rVert^{2} }
    \end{equation}
    
    As $\tilde{r_{k}}^{T}A^{-1}r_{k} - \tilde{r}_{k+1}^{T}A^{-1}r_{k+1} = \alpha_{k}\tilde{r}_{k}^{T}r_{k}$ for a non-symmetric matrix in BiCG:
    
    \begin{equation}
     \nonumber
        (T_{k+1}^{-1})_{(1,1)} = (T_{k}^{-1})_{(1,1)} + \frac{  ( \alpha_{k}\tilde{r}_{k}^{T}r_{k} ) } { \lVert r_{0} \rVert^{2} }
    \end{equation}
    
Thus knowing the relation between two consecutive $T_{k}$ inverse first elements we can relate it with the estimator we developed for A-norm of error and the approach of CGQL Algorithm. However relating Non-Symmetric Lanczos and BiCG through Quadrature based methods will involve Complex Gaussian Quadratures \cite{saylor2001gaussian}. Hence we have followed an equivalent but direct approach of estimation using the relations of Bi-CG explicitly.
    
\bibliographystyle{spmpsci}      % mathematics and physical sciences
%\bibliographystyle{siam}
\bibliography{reference}
\end{document}